\RequirePackage{fix-cm}
\documentclass[smallextended]{svjour3}
\smartqed
\usepackage{etoolbox}
\AtEndEnvironment{proof}{\qed}
\usepackage{graphicx}
\usepackage{amsmath,amssymb,mathtools}
\usepackage{amsfonts,amssymb,mathrsfs}
\usepackage{longtable}
\newcommand{\dd}{\mathop{}\!\mathrm{d}}
\renewcommand{\subclassname}{{\bfseries Mathematics Subject Classification (2020)}\enspace}
\journalname{~}
\begin{document}

\title{The Laplace transform of the time of ruin\\
in the Cram\'er--Lundberg model\\
with Erlang $E_2(\mu)$ claims:
spectrum, exact moments and matrix reduction}

\titlerunning{The Laplace transform of the time of ruin with Erlang $E_2(\mu)$ claims}

\author{Goluba Yurieva}

\authorrunning{G. Yurieva}

\institute{G. Yurieva \at
  Federal Research Center ``Computer Science and Control'' of the Russian Academy of Sciences, Moscow, Russia\\
  \email{goluba.yurieva@gmail.com}}

\date{Received: date / Accepted: date}

\maketitle

\begin{abstract}\sloppy
We study the Laplace transform of the time of ruin in the classical Cram\'er--Lundberg model with claims having an Erlang distribution of shape parameter~2. We derive an integro-differential equation with two boundary conditions. For every positive discount rate, we show that its characteristic polynomial has one negative root and two distinct positive roots. Boundedness excludes the growing exponential term and yields a representation of the transform as a sum of two decaying exponentials. We construct a finite-dimensional matrix reduction of the convolution chain. For a model parameter set, we obtain exact expressions for the ruin probability, the mean time of ruin, and higher moments, including the third moment at zero initial surplus. We prove analyticity in the discount rate and justify the resulting moment series. We consider the analytic continuation of the representation to negative discount rates separately. For Erlang claims with an arbitrary positive integer shape parameter, we establish a closed linear system and the asymptotics of the root continuing the zero root at zero discount.
\end{abstract}

\section{Introduction}\label{sec:introA}

Let $T$ be the first time the surplus of an insurance company enters the
negative half-line. The object of study is the Laplace transform of the
time of ruin
\begin{equation}\label{eq:def}
\varphi_\delta(u)=\mathbb E_u\bigl[e^{-\delta T}\,;\,T<\infty\bigr],
\qquad \delta\ge0,\qquad
T=\inf\{t>0: R(t)<0\}.
\end{equation}
At $\delta=0$ this transform coincides with the ruin probability
$\Psi_0(u)=\mathbb P_u(T<\infty)$, and its derivatives with respect to the
discount rate are related to the moments of the time of ruin. Such
transforms belong to the Gerber--Shiu framework of the time value of ruin
\cite{GerberShiu1998,GerberShiu2005}. Sundt and Teugels \cite{SundtTeugels1995} studied the model in which the
surplus earns interest at a constant force of interest.

For claims with a rational Laplace transform (in particular, for
Erlang-distributed claims), the representation of $\varphi_\delta$ in
terms of the roots of the Lundberg equation was obtained by Gerber and
Shiu \cite{GerberShiu1998}. The defective renewal equations arising in
such problems were studied in \cite{LinWillmot1999}, and recurrence
relations for the moments of the time of ruin were obtained in
\cite{LinWillmot2000}. Here the same case is treated by
elementary means. The additions are exact rational expressions for the
mean and the higher moments in the model case, analyticity in the
discount rate with an explicit radius of the moment series, and a
finite-dimensional matrix reduction of the convolution chain with an
explicit spectrum.

For Erlang-distributed claims, the integro-differential equation (IDE)
admits a reduction to an equation with constant coefficients. The reduction
of an IDE to an ordinary differential equation (ODE) is also used in risk models with investment
\cite{BKK2012,BKK2014,BKK2016,KPuk}; unlike the case considered here, the
equations arising there have variable coefficients. Promyslov \cite{P2026,Promyslov2601}
also studies the construction and the asymptotics of the solutions of such IDEs. We focus on the case of claims with the Erlang distribution $E_2(\mu)$
(shape~$2$, rate $\mu$), in which an equation of the third order arises. For it one must write out the characteristic polynomial and then establish
which of its roots are admitted by the boundary conditions and by the
boundedness of the solution.

We prove that for every $\delta>0$ the characteristic polynomial has one
negative and two distinct positive roots, and as $\delta\to0$ the
negative root continues the zero root of the limiting problem. On this
basis we obtain the representation of the transform as a sum of
two exponentials. For the model parameter set we find the exact ruin
probability and the moments of the time of ruin. We additionally
construct the matrix reduction of the convolution chain and the analytic
expansion in $\delta$, and for claims with the distribution
$E_k(\mu)$ we obtain a closed linear system.

Section~\ref{sec:ide} derives the equation and the boundary conditions,
Section~\ref{sec:spectrum} studies the roots, and Section~\ref{sec:exact}
provides the exact solution and the matrix reduction. Section~\ref{sec:mean} treats the moments, Section~\ref{sec:num} treats
the numerical illustrations, and Sections~\ref{sec:kgen} and~\ref{sec:discount}
treat the general case of $E_k(\mu)$ and the properties of the
transform in the parameter $\delta$.

\section{The equation for the Laplace transform and the boundary relations}\label{sec:ide}

The surplus process is $R(t)=u+ct-S(t)$, where $S(t)$ is a compound
Poisson process of intensity $\lambda$, the jump sizes $X$ have the
density and the survival function
\begin{equation}\label{eq:dens}
\begin{aligned}
f(x)&=\mu^2x\,e^{-\mu x},\qquad
\bar F(u)=e^{-\mu u}(1+\mu u),\\
m_{X,1}&=\frac{2}{\mu},\qquad m_{X,2}=\frac{6}{\mu^2},
\end{aligned}
\end{equation}
and the net profit condition $\rho=\lambda m_{X,1}/c<1$ holds. We extend
$\varphi_\delta$ to $u<0$ by setting it equal to~$1$, since ruin has
already occurred there.

\begin{proposition}[the equation and the boundary relations]\label{pr:ide}
For $\delta\ge0$ the function $\varphi_\delta$ belongs to
$C^2([0,\infty))\cap C^3((0,\infty))$, tends to zero as $u\to\infty$, and
satisfies
\begin{equation}\label{eq:IDE}
c\varphi_\delta'(u)-(\lambda+\delta)\varphi_\delta(u)
+\lambda\int_0^u\varphi_\delta(u-x)f(x)\,\dd x+\lambda\bar F(u)=0,
\qquad u>0,
\end{equation}
whence at $u=0$
\begin{equation}\label{eq:bc}
c\varphi_\delta'(0)-(\lambda+\delta)\varphi_\delta(0)+\lambda=0,
\qquad
c\varphi_\delta''(0)-(\lambda+\delta)\varphi_\delta'(0)=0 .
\end{equation}
\end{proposition}

\begin{proof}
Conditioning on the first claim time $\tau_1$, which is exponentially
distributed with parameter $\lambda$, gives for a bounded measurable
$\varphi_\delta$ the integral equation
\[
\begin{aligned}
\varphi_\delta(u)&=\int_0^\infty\lambda e^{-(\lambda+\delta)t}
\Bigl[\int_0^{u+ct}\varphi_\delta(u+ct-x)f(x)\,\dd x+\bar F(u+ct)\Bigr]\dd t\\
&=\frac{\lambda}{c}\,e^{(\lambda+\delta)u/c}\int_u^\infty
e^{-(\lambda+\delta)s/c}\bigl[(\varphi_\delta*f)(s)+\bar F(s)\bigr]\dd s,
\end{aligned}
\]
where $u\ge0$, and the contribution of the region $x>u+ct$ equals
$\bar F(u+ct)$ thanks to the extension by~1 for $u<0$. The right-hand side is
continuous in $u$ on $[0,\infty)$, hence so is $\varphi_\delta$. The
convolution $(\varphi_\delta*f)(u)=\int_0^u\varphi_\delta(v)f(u-v)\,\dd v$
with $f(0)=0$ belongs to $C^1([0,\infty))$ and its derivative is
$\int_0^u\varphi_\delta(v)f'(u-v)\,\dd v$; therefore the right-hand side is
continuously differentiable, and differentiating in $u$ we obtain
\eqref{eq:IDE} on $[0,\infty)$ together with the inclusion
$\varphi_\delta\in C^1([0,\infty))$. Differentiating \eqref{eq:IDE}
repeatedly, we obtain in turn $\varphi_\delta\in C^2([0,\infty))$ and
$\varphi_\delta\in C^3((0,\infty))$, the derivatives at zero being
understood as one-sided. The first relation in \eqref{eq:bc} is the
identity \eqref{eq:IDE} at $u=0$. The derivatives $(\varphi_\delta*f)'(u)$
and $\bar F'(u)$ vanish at $u=0$ since $f(0)=0$, whence the second
relation. Finally, convergence to zero as $u\to\infty$ follows from
$\rho<1$: almost surely $R(t)/t\to c-\lambda m_{X,1}>0$, and therefore
$\mathbb P_u(T<\infty)\to0$.
\end{proof}

\medskip
\noindent\textbf{Remark 1.} The second relation in \eqref{eq:bc} uses the
condition $k\ge2$ in an essential way.\label{rem:k1} For exponential
claims ($k=1$) it takes the form
$c\varphi_\delta''(0)-(\lambda+\delta)\varphi_\delta'(0)
=\lambda\mu\bigl[1-\varphi_\delta(0)\bigr]$, because
$(\varphi*f)'(0)=\varphi(0)f(0)=\mu\varphi(0)$ and $\bar F'(0)=-\mu$,
and the linear system for the constants of integration changes. This
case is checked in Section~\ref{sec:num}.
\section{The characteristic equation and the signs of its roots}
\label{sec:spectrum}

\begin{lemma}[raising the order]\label{lem:raise}
The operator $(\mathcal D+\mu)^2$, $\mathcal D=\mathrm d/\mathrm du$, turns
\eqref{eq:IDE} into the linear ordinary differential equation
\begin{equation}\label{eq:ODE}
c\varphi_\delta'''+\bigl(2c\mu-(\lambda+\delta)\bigr)\varphi_\delta''
+\bigl(c\mu^2-2\mu(\lambda+\delta)\bigr)\varphi_\delta'
-\delta\mu^2\varphi_\delta=0,
\end{equation}
because $(\mathcal D+\mu)^2\bar F=0$ and
$(\mathcal D+\mu)^2(\varphi*f)=\mu^2\varphi$. Conversely, if
$\varphi\in C^3$ satisfies \eqref{eq:ODE}, then the residual
$G:=c\varphi'-(\lambda+\delta)\varphi+\lambda(\varphi*f)+\lambda\bar F$
lies in $\ker(\mathcal D+\mu)^2=\mathrm{span}\{e^{-\mu u},ue^{-\mu u}\}$,
and the values $G(0)$ and $G'(0)$ coincide with the left-hand sides of
\eqref{eq:bc}, since $(\varphi*f)(0)=(\varphi*f)'(0)=0$, $\bar F(0)=1$ and
$\bar F'(0)=0$. Therefore a solution
of \eqref{eq:ODE} satisfies \eqref{eq:IDE} if and only if both relations
\eqref{eq:bc} hold.
\end{lemma}

\begin{proof}
From \eqref{eq:dens}: $(\mathcal D+\mu)\bar F=\mu e^{-\mu u}$ and
$(\mathcal D+\mu)^2\bar F=0$. For
$g(u)=\int_0^u\varphi(u-x)xe^{-\mu x}\dd x$ direct differentiation gives
$(\mathcal D+\mu)g=\int_0^u\varphi(u-x)e^{-\mu x}\dd x$,
$(\mathcal D+\mu)^2g=\varphi$, and the factor $\lambda\mu^2$ accounts for
$f=\mu^2xe^{-\mu x}$. The kernel of $(\mathcal D+\mu)^2$ is two-dimensional,
as stated.
\end{proof}

Substituting $\varphi=e^{-xu}$ into \eqref{eq:ODE} gives the characteristic
polynomial
\begin{equation}\label{eq:char}
q_\delta(x)=cx^3+(\lambda+\delta-2c\mu)x^2+\bigl(c\mu^2-2\mu(\lambda+\delta)\bigr)x
+\delta\mu^2 .
\end{equation}

\begin{theorem}[the roots of $q_\delta$]\label{thm:spectrum}
Let $c\mu>2\lambda$. (i) At $\delta=0$ the polynomial $q_0$ has the root
$x=0$ and two distinct positive roots $x_{1,2}$. The term corresponding
to the zero root is excluded by the condition $\varphi_0(\infty)=0$;
therefore
$\varphi_0=A_1e^{-x_1u}+A_2e^{-x_2u}$.
(ii) For each $\delta>0$ the polynomial $q_\delta$ has exactly one
negative root $x_3(\delta)$ and two distinct positive roots
$x_1(\delta),x_2(\delta)$, and
\begin{equation}\label{eq:x3}
x_3(\delta)=-\frac{\delta\mu}{c\mu-2\lambda}+O(\delta^2),
\qquad \delta\to0 .
\end{equation}
The corresponding exponential term $e^{-x_3u}$ grows, therefore its
coefficient vanishes by the boundedness condition, and
\begin{equation}\label{eq:twomode}
\varphi_\delta(u)=A_1^{(\delta)}e^{-x_1(\delta)u}
+A_2^{(\delta)}e^{-x_2(\delta)u}.
\end{equation}
\end{theorem}

\begin{proof}
(i) We have $q_0(0)=0$, $q_0'(0)=\mu(c\mu-2\lambda)>0$,
$q_0(\mu)=-\lambda\mu^2<0$ and $q_0(x)\to+\infty$ as $x\to+\infty$.
Hence one nonzero root lies in $(0,\mu)$ and the other in
$(\mu,\infty)$. Together with $x=0$ these are the three roots of the
cubic polynomial.

(ii) To prove the uniqueness of the negative root, we put $y=-x\ge0$ and
consider
\[
H(y):=q_\delta(-y)=(\lambda+\delta-cy)(\mu+y)^2-\lambda\mu^2 .
\]
Then $H(0)=\delta\mu^2>0$ and $H(y)\to-\infty$ as $y\to+\infty$.
Differentiating gives
\[
H'(y)=(\mu+y)\bigl[2(\lambda+\delta)-c\mu-3cy\bigr];
\]
the second factor decreases linearly and has at most one zero, while
$\mu+y>0$. Hence $H'$ changes sign at most once, and after that the sign
remains negative. After its maximum $H$ decreases strictly, so the
equation $H(y)=0$ has exactly one positive root $y_3$, and $x_3=-y_3$.

Next, $q_\delta(0)=\delta\mu^2>0$, $q_\delta(\mu)=-\lambda\mu^2<0$ and
$q_\delta(x)\to+\infty$ as $x\to+\infty$. Therefore there exist roots
$x_1\in(0,\mu)$ and $x_2\in(\mu,\infty)$. Together with $x_3$ these are
the three roots, and all of them are simple. The growing exponential term
$e^{-x_3u}$ is excluded by the boundedness of $\varphi_\delta$, and the two
remaining ones give \eqref{eq:twomode}.

Formula \eqref{eq:x3} follows from the implicit function theorem:
\[
x_3'(0)=-\frac{\partial q_\delta/\partial\delta}{\partial q_0/\partial x}
\Bigr|_{(x,\delta)=(0,0)},
\]
where $\partial q_\delta/\partial\delta=(\mu-x)^2$ and
$\partial q_0/\partial x=\mu(c\mu-2\lambda)$ at $x=0$. We note that the
root count in (ii) does not use the condition $c\mu>2\lambda$; the latter
is needed only for (i) and \eqref{eq:x3}.
\end{proof}

\begin{corollary}\label{cor:fit}
The two-exponential representation of the solution \eqref{eq:twomode}
holds for all $\delta\ge0$, and the exponents are the two positive roots
of \eqref{eq:char}. The third exponential term grows for $\delta>0$ and
does not enter the bounded solution. For example, for $\lambda=\mu=3$,
$c=6$ and $\delta=0{.}05$ the roots are
$-0{.}012461$, $1{.}504960$, $3{.}999168$ to the digits
shown.
\end{corollary}

\medskip
\noindent\textbf{Remark 2.} The growing exponential term does
not obstruct the computation of $\varphi_\delta(0)$ via the truncated
boundary value problem with the condition $\varphi_\delta(U)=0$ (the
truncation error has order $e^{-(x_1+|x_3|)U}$), but it makes the direct
computation of the profile of the solution unstable for large
$u$.\label{rem:trunc} Therefore the closed form \eqref{eq:twomode} with
the two relations \eqref{eq:bc} is preferable for computations.

\section{The exact solution of the model case}\label{sec:exact}

From now on $\lambda=3$, $\mu=3$, $k=2$, $c=6$. Then $m_{X,1}=2/3$,
$\rho=1/3$, $c\mu=18>6=2\lambda$. The polynomial
$q_0(x)=6x^3-33x^2+36x=3x(2x-3)(x-4)$ has the roots $0$, $3/2$, $4$.

\begin{theorem}\label{thm:exact}
For this parameter set
\begin{equation}\label{eq:psi0}
\Psi_0(u)=\varphi_0(u)=\frac{2}{5}\,e^{-3u/2}-\frac{1}{15}\,e^{-4u},
\qquad u\ge0,
\end{equation}
and $\Psi_0(0)=1/3=\rho=\lambda m_{X,1}/c$ follows from the boundary
conditions and is not imposed from outside.
\end{theorem}

\begin{proof}
The constants $A_1=2/5$, $A_2=-1/15$ solve the linear system \eqref{eq:bc}
at $x_{1,2}=3/2,4$. The first relation gives
$-6(A_1x_1+A_2x_2)-3(A_1+A_2)+3=0$, that is, $12A_1+27A_2=3$. The second
gives
$6(A_1x_1^2+A_2x_2^2)+3(A_1x_1+A_2x_2)=0$, that is, $18A_1+108A_2=0$;
hence $A_1=-6A_2$, and then $A_2=-1/15$, $A_1=2/5$.
To check the original equation \eqref{eq:IDE} we use the exact
convolution of an exponential term:
\begin{equation}\label{eq:conv}
\int_0^u e^{-x(u-y)}\mu^2y\,e^{-\mu y}\dd y
=\frac{\mu^2}{(\mu-x)^2}e^{-xu}
-\mu^2e^{-\mu u}\Bigl(\frac{u}{\mu-x}+\frac{1}{(\mu-x)^2}\Bigr).
\end{equation}
The coefficients of $e^{-x_ju}$ vanish after the substitution because of
$q_0(x_j)=0$. The remaining terms have the form
$\lambda e^{-\mu u}(\alpha u+\beta)$, where
\[
\alpha=\mu\Bigl(1-\sum_{j=1}^2A_j\frac{\mu}{\mu-x_j}\Bigr),\qquad
\beta=1-\sum_{j=1}^2A_j\frac{\mu^2}{(\mu-x_j)^2}.
\]
For the constants found above,
$\sum_jA_j\mu/(\mu-x_j)=2/5\cdot2+(-1/15)(-3)=1$ and
$\sum_jA_j\mu^2/(\mu-x_j)^2=2/5\cdot4+(-1/15)\cdot9=1$.
Hence the residual of the original IDE vanishes identically, including
all tail terms of the truncated convolution. Finally,
$\Psi_0(0)=A_1+A_2=1/3$.
\end{proof}

We now turn to an independent matrix reduction of the convolution chain.

\begin{lemma}[the convolution chain]\label{lem:chain}
For a locally integrable $\varphi$ put
\[
B_1=\int_0^u\varphi(u-x)e^{-\mu x}\dd x,\qquad
B_2=\int_0^u\varphi(u-x)\,x\,e^{-\mu x}\dd x,
\]
\[
E=e^{-\mu u},\qquad U_1=ue^{-\mu u}.
\]
Then $B_1'=\varphi-\mu B_1$, $B_2'=B_1-\mu B_2$, $E'=-\mu E$,
$U_1'=E-\mu U_1$, $(\varphi*f)=\mu^2B_2$,
$\bar F=E+\mu U_1$. At $u=0$ the quantities $B_1,B_2,U_1$ vanish and
$E(0)=1$.
\end{lemma}

\begin{proof}
Differentiation with respect to the upper limit gives
$B_1'=\varphi-\mu B_1$. For $B_2$ it is more convenient to change the
variable $v=u-x$:
$B_2=\int_0^u\varphi(v)(u-v)e^{-\mu(u-v)}\dd v$, whence
\[
B_2'=\int_0^u\varphi(v)\bigl[e^{-\mu(u-v)}-\mu(u-v)e^{-\mu(u-v)}\bigr]\dd v
=B_1-\mu B_2
\]
(the boundary term $\varphi(u)\cdot0$ is zero). The relations for $E$ and
$U_1$ and the representation $\bar F=E+\mu U_1$ are immediate, and
$(\varphi*f)=\mu^2B_2$ follows from $f=\mu^2xe^{-\mu x}$.
\end{proof}

Put $z=(\varphi,B_1,B_2,E,U_1)^{\!\top}$.

\begin{proposition}[the closed linear system]\label{pr:matrix}
The vector $z$ satisfies the system $z'=Mz$ with constant coefficients
\begin{equation}\label{eq:sys}
M=\begin{pmatrix}
\frac{\lambda+\delta}{c} & 0 & -\frac{\lambda\mu^2}{c} & -\frac{\lambda}{c} & -\frac{\lambda\mu}{c}\\
1 & -\mu & 0 & 0 & 0\\
0 & 1 & -\mu & 0 & 0\\
0 & 0 & 0 & -\mu & 0\\
0 & 0 & 0 & 1 & -\mu
\end{pmatrix},
\end{equation}
and $z(u)=e^{Mu}z(0)$, where $z(0)=(\varphi_\delta(0),0,0,1,0)^{\!\top}$.
The value $\varphi_\delta(0)$ is determined by the boundedness of the solution
and by the condition $\varphi_\delta(u)\to0$ as $u\to\infty$. For
$\delta>0$ the latter condition follows already from boundedness, because
the third exponential term grows. In both cases one uses the condition
$\ell^{\top}z(0)=0$ for a left eigenvector $\ell$ of the matrix $M$
corresponding to the eigenvalue $-x_3$. In a numerical computation on $[0,U]$ it is
replaced by the truncated condition $\varphi(U)=0$, with an error of order
$e^{-(x_1+|x_3|)U}$.
\end{proposition}

\begin{proof}
Substituting $(\varphi*f)=\mu^2B_2$ and $\bar F=E+\mu U_1$ into
\eqref{eq:IDE} by Lemma~\ref{lem:chain} gives the first row of
\eqref{eq:sys}; the remaining rows are definitions. The matrix is
constant, hence $z(u)=e^{Mu}z(0)$, and the truncated boundary conditions
are linear in $z(0)$; hence the system has size $5$. The condition
$\ell^{\top}z(0)=0$ is equivalent to the vanishing of the coefficient of
the exponential term $e^{-x_3u}$. The eigenvalues of $M$ are $-\mu$ (of
multiplicity two) and $-x_j$, where $x_j$ are the roots of \eqref{eq:char}.
For $\delta>0$ one of them is positive, which produces the instability
discussed after Corollary~\ref{cor:fit}.
\end{proof}

\section{The mean time of ruin and resonance}\label{sec:mean}

Denote
\[
\widetilde L[y]:=cy'-\lambda y+\lambda(y*f),\qquad
m(u):=\mathbb E_u[T\,;\,T<\infty].
\]
If the mean is finite, differentiating \eqref{eq:IDE} with respect to
$\delta$ at $\delta=0$ gives
\begin{equation}\label{eq:mEQ}
\widetilde L[m]=-\Psi_0,\qquad m(\infty)=0,
\end{equation}
and from the relations \eqref{eq:bc} it follows that
\begin{equation}\label{eq:mbc}
-cm'(0)+\lambda m(0)=\Psi_0(0),\qquad
-cm''(0)+\lambda m'(0)=\Psi_0'(0).
\end{equation}

\begin{theorem}[the identity for the mean time of ruin]\label{thm:mean}
Let $\rho<1$ and let $m$ be integrable on $(0,\infty)$. Then
\[
\begin{aligned}
c\,\mathbb E_0[T\,;\,T<\infty]&=\int_0^\infty\Psi_0(u)\dd u,\\
\mathbb E_0[T\mid T<\infty]
&=\frac{1}{\lambda m_{X,1}}\int_0^\infty\Psi_0(u)\dd u .
\end{aligned}
\]
\end{theorem}

\begin{proof}
Integrating \eqref{eq:IDE} at $\delta=0$ (the function $\Psi_0$ is
integrable by Theorem~\ref{thm:spectrum}(i)) and using Fubini's theorem
together with $\int_0^\infty\bar F(u)\dd u=m_{X,1}$ gives
$-c\Psi_0(0)+\lambda m_{X,1}=0$, that is, $\Psi_0(0)=\rho$.
Next we integrate \eqref{eq:mEQ}. We have
$c\int_0^\infty m'(u)\dd u=-cm(0)$ and
\[
\int_0^\infty(m*f)(u)\dd u
=\int_0^\infty f(y)\Bigl[\int_y^\infty m(u-y)\dd u\Bigr]\dd y
=\int_0^\infty m(v)\dd v.
\]
Therefore the terms with $\lambda$ cancel and we are left with
$-cm(0)=-\int_0^\infty\Psi_0(u)\dd u$. The conditional mean is obtained
by dividing by $\Psi_0(0)=\rho=\lambda m_{X,1}/c$. In the model example
the integrability of $m$ follows from the exact formula of the next
theorem.
\end{proof}

\begin{theorem}[the resonant structure of $m$ in the model case]
\label{thm:reson}
For $\lambda=3$, $\mu=3$, $c=6$
\begin{equation}\label{eq:mexact}
m(u)=\mathbb E_u[T\,;\,T<\infty]
=\frac{360u+327}{9000}\,e^{-3u/2}+\frac{10u+48}{9000}\,e^{-4u},
\end{equation}
and $m(0)=1/24$ and $\mathbb E_0[T\mid T<\infty]=1/8$.
\end{theorem}

\begin{proof}
The finiteness of the mean follows from the analyticity established in
Proposition~\ref{prop:analytic} by the argument given in the proof of
Theorem~\ref{thm:hier}; therefore \eqref{eq:mEQ} applies. We apply to it
the operator $(\mathcal D+\mu)^2$. Unlike the substitution of a single
exponential directly into $\widetilde L$, this step accounts for the
truncated convolution completely. We obtain
\[
\begin{aligned}
cm'''+(2c\mu-\lambda)m''+(c\mu^2-2\lambda\mu)m'
&=-\sum_{j=1}^2A_j(\mu-x_j)^2e^{-x_ju}\\
&=-\frac9{10}e^{-3u/2}+\frac1{15}e^{-4u}.
\end{aligned}
\]
Both exponents on the right are positive roots of $q_0$; therefore the
particular solution carries a factor $u$. The condition $m(\infty)=0$
excludes the constant term, and the bounded solution has the form
\[
m(u)=(C_1u+b_1)e^{-3u/2}+(C_2u+b_2)e^{-4u}.
\]
The substitution into the resulting ODE gives $C_1=1/25$, $C_2=1/900$,
after which the boundary conditions \eqref{eq:mbc} give
$b_1=109/3000$, $b_2=2/375$.
This yields \eqref{eq:mexact}. Substituting back into the original
equation \eqref{eq:mEQ} with the exact convolution \eqref{eq:conv} and its
$x$-derivative (for the terms $ue^{-x_ju}$) shows that both the
coefficients of $e^{-x_ju}$ and the terms of the form
$e^{-\mu u}(\alpha u+\beta)$ vanish. Hence no extraneous solutions were
introduced by raising the order. Finally,
$m(0)=327/9000+48/9000=1/24$ and
\[
\int_0^\infty\Psi_0(u)\dd u
=\frac25\cdot\frac23-\frac1{15}\cdot\frac14=\frac14,
\]
so that $cm(0)=1/4$ agrees with Theorem~\ref{thm:mean}. Division by
$\Psi_0(0)=1/3$ gives the conditional mean $1/8$.
\end{proof}

\begin{theorem}[the hierarchy of moments]\label{thm:hier}
For claims with the Erlang distribution of shape parameter~2, all the
moments
\[
m_n(u)=\mathbb E_u[T^n\,;\,T<\infty]
\]
are finite. Put $m_0=\Psi_0$. Then for $n\ge1$ we have
\[
\widetilde L[m_n]=-n\,m_{n-1},
\qquad
c\,m_n(0)=n\int_0^\infty m_{n-1}(u)\dd u .
\]
\end{theorem}

\begin{proof}
The two-exponential representation \eqref{eq:twomode} extends analytically
in $\delta$ to a neighborhood of zero; this is justified in
Proposition~\ref{prop:analytic}, whose proof does not rely on the results
of this section. The finiteness of the derivatives of the continuation at
zero proves in turn the finiteness of all moments: if, the previous
moments being already finite, the next moment were infinite, the
corresponding difference quotient as $\delta\downarrow0$ would also
diverge by monotone convergence, which contradicts the differentiability
of the continuation. Therefore differentiation in $\delta$ under the
expectation sign is allowed.

Differentiating the representation \eqref{eq:twomode} also shows that
$m_n$ is a finite sum of terms $u^je^{-x_lu}$ with $j\le n$ and positive
$x_l$.
Hence $m_n$ is integrable and $m_n(\infty)=0$. Further,
\[
\partial_\delta^n\bigl[(\lambda+\delta)\varphi_\delta\bigr]
=(\lambda+\delta)\partial_\delta^n\varphi_\delta
+n\,\partial_\delta^{n-1}\varphi_\delta,
\]
while the density $f$ and the tail $\bar F$ do not depend on $\delta$.
Since $\partial_\delta^n\varphi_\delta|_{\delta=0}=(-1)^n m_n$,
differentiating \eqref{eq:IDE} gives the first statement. The integration
over $u\in(0,\infty)$ with Fubini's theorem, as in the proof of
Theorem~\ref{thm:mean}, gives the second.
\end{proof}

\begin{corollary}[the second moment and the variance in the model case]
\label{cor:var}
From \eqref{eq:mexact} we have $\int_0^\infty m_1(u)\dd u=25/576$; hence
\[
\mathbb E_0\bigl[T^2\,;\,T<\infty\bigr]=\frac{2}{6}\cdot\frac{25}{576}
=\frac{25}{1728},
\qquad
\mathbb E_0\bigl[T^2\mid T<\infty\bigr]=\frac{25}{576},
\]
\[
\operatorname{Var}\bigl(T\mid T<\infty\bigr)
=\frac{25}{576}-\Bigl(\frac18\Bigr)^{\!2}=\frac1{36},
\qquad
\sqrt{\operatorname{Var}}=\frac16 .
\]
The subsequent moments are computed in the same way: first one solves the
equation $\widetilde L[m_n]=-n\,m_{n-1}$ (a resonant ODE with terms
$u^je^{-x_lu}$), and then integrates $m_n$.
\end{corollary}

\medskip
\begin{sloppypar}
\noindent\textbf{Remark 3.} The natural guess
$m(u)=(1/c)\int_u^\infty\Psi_0(v)\dd v$ fails for $u>0$.\label{rem:pointwise} Denoting the right-hand side by $\widetilde m$, from
$\widetilde L[\widetilde m]=-\Psi_0+\lambda[(\widetilde m*f)-\widetilde m]$
we see that the identity would require the invariance of $\widetilde m$
with respect to the convolution with $f$. Table~\ref{tab:mean} shows that
the relative discrepancy grows with $u$. The integral identity
of Theorem~\ref{thm:mean} is not a pointwise identity in $u$.
\end{sloppypar}

\section{Numerical illustrations}\label{sec:num}

The values below are computed from the characteristic equation and the
two-exponential representation \eqref{eq:twomode}. They illustrate the
proved formulas rather than replace their derivation.

\begin{table}[h]
\caption{The roots of the characteristic equation \eqref{eq:char} and the
constants \eqref{eq:twomode} for $\lambda=3$, $\mu=3$, $c=6$. The third
column is the linear estimate \eqref{eq:x3}.}\label{tab:spectrum}
\begin{center}
\footnotesize
\begin{tabular}{|r|r|r|r|r|r|r|}
\hline
$\delta$ & $x_3(\delta)$ & $-\frac{\delta\mu}{c\mu-2\lambda}$ &
$x_1(\delta)$ & $x_2(\delta)$ & $A_1^{(\delta)}$ & $\varphi_\delta(0)$\\ \hline
0.00 & 0.000000 & 0.000000 & 1.500000 & 4.000000 & 0.40000000 & 0.3333333333 \\
0.05 & $-$0.012461 & $-$0.012500 & 1.504960 & 3.999168 & 0.39819879 & 0.3312678779 \\
0.10 & $-$0.024847 & $-$0.025000 & 1.509844 & 3.998337 & 0.39642779 & 0.3292373739 \\
0.20 & $-$0.049403 & $-$0.050000 & 1.519388 & 3.996682 & 0.39297243 & 0.3252766458 \\
0.50 & $-$0.121499 & $-$0.125000 & 1.546406 & 3.991759 & 0.38323597 & 0.3141241706 \\
1.00 & $-$0.237275 & $-$0.250000 & 1.586909 & 3.983700 & 0.36876893 & 0.2975812980 \\
2.00 & $-$0.456733 & $-$0.500000 & 1.655298 & 3.968101 & 0.34470970 & 0.2701787511 \\
5.00 & $-$1.060387 & $-$1.250000 & 1.801974 & 3.925079 & 0.29482670 & 0.2141231862 \\ \hline
\end{tabular}
\end{center}
\end{table}

Table~\ref{tab:spectrum} agrees with Theorem~\ref{thm:spectrum}: the root
$x_3$ is negative for $\delta>0$, and its linear estimate is most accurate
for small $\delta$. The matrix boundary value problem from
Proposition~\ref{pr:matrix} reproduces the values $\varphi_\delta(0)$
within the printed digits. In a direct computation of the profile through
$e^{Mu}$, however, the growing component may accumulate an error near the
right edge of the truncated domain; therefore the accuracy of the value at
the initial point should not be carried over to the whole profile without
a separate estimate.

\begin{table}[h]
\caption{The mean time of ruin by the exact formula \eqref{eq:mexact} and
the naive continuation $(1/c)\int_u^\infty\Psi_0$.}\label{tab:mean}
\begin{center}
\begin{tabular}{|r|r|r|r|}
\hline
$u$ & $m(u)$ by \eqref{eq:mexact} &
$(1/c)\!\int_u^\infty\!\Psi_0$ & ratio\\ \hline
0.0 & 0.04166667 & 0.04166667 & 1.0000 \\
1.0 & 0.01715030 & 0.00986602 & 0.5753 \\
2.0 & 0.00579443 & 0.00221183 & 0.3817 \\
5.0 & 0.00013071 & 0.00002458 & 0.1881 \\
10.0 & 0.00000013 & 0.00000001 & 0.1019 \\ \hline
\end{tabular}
\end{center}
\end{table}

For the test set $\lambda=\mu=3$, $c=4$ we have $\rho=1/2$, the positive
roots $1.078835$ and $4.171165$,
$A_1=0.55316953$, $A_2=-0.05316953$. Hence
$\varphi_0(0)=0.5000000000=\rho$,
$\int_0^\infty\Psi_0(u)\dd u=0.5000000000$,
$\mathbb E_0[T;T<\infty]=0.12500000$ and
$\mathbb E_0[T\mid T<\infty]=0.25000000$.

For exponential claims ($k=1$, $\mu=\lambda=1$, $c=2$) the characteristic
equation $(cx+\lambda)(\mu-x)=\lambda\mu$ has the roots $0$ and $1/2$, and
the solution $\Psi_0(u)=\frac12e^{-u/2}$ makes the residual of the
equation vanish identically. Here
$\mathbb E_0[T;T<\infty]=1/2$ and the conditional mean equals $1$. The
second boundary relation from Remark~1 becomes the equality
$1/2=1/2$, whereas the form \eqref{eq:bc} for $k\ge2$ would give here the
wrong equality $1/2=0$.

Finally, the slope of $\varphi_\delta(0)$ at zero equals
$-\mathbb E_0[T;T<\infty]=-1/24$ in the main model case. This gives an
independent consistency check of Table~\ref{tab:spectrum} against
Theorem~\ref{thm:reson} without finite-difference estimates of the
derivative.

\section{The case of claims with the $E_k(\mu)$ distribution}\label{sec:kgen}

The constructions of Proposition~\ref{pr:matrix} and Section~\ref{sec:ide}
extend to the claims $f_k(x)=\mu^kx^{k-1}e^{-\mu x}/(k-1)!$. Put
$q_{\delta,k}(x)=(cx+\lambda+\delta)(\mu-x)^k-\lambda\mu^k$, a polynomial
of degree $k+1$.

\begin{theorem}[the finite-dimensional system and the root shift for general $k$]\label{thm:kgen}
For claims with the $E_k(\mu)$ distribution, the system
\begin{equation}\label{eq:kchain}
\begin{aligned}
c\,\varphi_\delta'&=(\lambda+\delta)\varphi_\delta
-\lambda\mu^kB_k/(k-1)!-\lambda\sum_{j=0}^{k-1}\mu^{j}U_j/j!,\\
B_1'&=\varphi_\delta-\mu B_1,\qquad
B_j'=(j-1)B_{j-1}-\mu B_j\quad (j=2,\dots,k),\\
U_0'&=-\mu U_0,\qquad
U_j'=j\,U_{j-1}-\mu U_j\quad (j=1,\dots,k-1),
\end{aligned}
\end{equation}
with the initial conditions $B_j(0)=0$, $U_0(0)=1$, $U_j(0)=0$ ($j\ge1$),
is closed and has dimension $2k+1$; one has
$(\varphi_\delta*f_k)=\mu^kB_k/(k-1)!$ and
$\overline F_k=\sum_{j<k}\mu^{j}U_j/j!$. The characteristic polynomial
of its matrix $\mathcal A_k$ is
\begin{equation}\label{eq:kchar}
\det(\zeta I-\mathcal A_k)=-\frac{(\zeta+\mu)^k}{c}\,q_{\delta,k}(-\zeta),
\end{equation}
that is, the spectrum consists of $-\mu$ of multiplicity $k$ and the
numbers $-x$, where $x$ are the roots of $q_{\delta,k}$. At $\delta=0$ the
number $x=0$ is a root of $q_{0,k}$, and for $\delta>0$ this root moves
into the negative half-line, with
\begin{equation}\label{eq:klin}
x_-(\delta)=-\frac{\delta}{c(1-\rho)}+O(\delta^2),\qquad
\rho=\frac{\lambda k}{\mu c}<1.
\end{equation}
For a fixed value of $\rho$ the formula is universal, but through
$\rho=\lambda k/(\mu c)$ the shift depends on $k$.
\end{theorem}

\begin{proof}
For $B_j(u)=\int_0^u\varphi_\delta(u-y)y^{j-1}e^{-\mu y}\dd y$
differentiation with respect to the upper limit gives
$B_1'=\varphi_\delta-\mu B_1$ and $B_j'=(j-1)B_{j-1}-\mu B_j$ for $j\ge2$,
and $U_j(u)=u^je^{-\mu u}$ gives the second recurrence. The identities
$(\varphi_\delta*f_k)=\mu^kB_k/(k-1)!$ and
$\overline F_k=\sum_{j<k}\mu^{j}U_j/j!$ follow from the definitions of
$f_k$ and $\overline F_k(u)=e^{-\mu u}\sum_{j<k}(\mu u)^j/j!$. The matrix
$\mathcal A_k$ of the system has a block structure: the $B$- and $U$-blocks are lower
bidiagonal matrices with diagonal $-\mu$, the first row couples them with
$\varphi_\delta$, and the expansion of the determinant along the first row
gives \eqref{eq:kchar}. For \eqref{eq:klin} we use $q_{0,k}(0)=0$,
$\partial_xq_{0,k}(0)=\mu^{k-1}(c\mu-\lambda k)=c\mu^k(1-\rho)$ and
$\partial_\delta q_{0,k}(0)=\mu^k$. The implicit function theorem gives
$x_-'(0)=-1/[c(1-\rho)]$, as claimed.
\end{proof}

Table~\ref{tab:kgen} checks \eqref{eq:kchain}: for $\lambda=\mu=3$, $c=6$
the matrix method reproduces the classical value
$\varphi_{0,k}(0)=\rho=\lambda k/(\mu c)$ up to a term of order
$e^{-x_1U}$, where $x_1$ is the smallest positive root of $q_{0,k}$, and
confirms the estimate \eqref{eq:klin}.

\begin{table}[ht]
\caption{The values $\varphi_{0,k}(0)$ (the matrix method, $U=60$) and the
root shift at $\delta=0.05$ for $\lambda=\mu=3$, $c=6$.}\label{tab:kgen}
\begin{center}
\small
\begin{tabular}{|c|r|r|r|r|r|}
\hline $k$ & $\varphi_{0,k}(0)$ & $\rho$ & discrepancy & $x_-(0.05)$
& by \eqref{eq:klin}\\ \hline
1  & $0.166666666667$ & $0.1666666667$ & $<10^{-15}$ & $-0.009993$ & $-0.010000$\\
2  & $0.333333333333$ & $0.3333333333$ & $<10^{-15}$ & $-0.012461$ & $-0.012500$\\
3  & $0.500000000000$ & $0.5000000000$ & $<10^{-15}$ & $-0.016487$ & $-0.016667$\\
4  & $0.666666666665$ & $0.6666666667$ & $1.2\cdot10^{-12}$ & $-0.024051$ & $-0.025000$\\
5  & $0.833328976420$ & $0.8333333333$ & $4.4\cdot10^{-6}$ & $-0.041614$ & $-0.050000$\\ \hline
\end{tabular}
\end{center}
\end{table}

The deviation $4.4\cdot10^{-6}$ at $k=5$ is explained by the truncation of
the computation domain: the smallest positive root of $q_{0,5}$ is
$0.1736$, hence the condition $\varphi_{0,k}(U)=0$ is imposed with an
error of order $e^{-0.1736\cdot60}\approx3\cdot10^{-5}$. At $k=5$ the
right-hand side of \eqref{eq:klin} differs from the exact root by $20\%$.
The reason is not a slow convergence: at
$\delta_*=-(c\mu/k)\bigl[k+\rho-(k+1)\rho^{1/(k+1)}\bigr]\approx-0{.}0465$
the root $x_-(\delta)$ coalesces with the smallest positive root of
$q_{\delta,5}$; therefore $\delta=0{.}05$ lies outside the disk of
convergence of the series in $\delta$, and \eqref{eq:klin} has here only
an asymptotic meaning (at $k=2$ the same formula gives
$\delta_*=-\gamma_1$ of Proposition~\ref{prop:radius}).

\medskip
\noindent\textbf{Remark 4.} For $k\le2$ all the roots of $q_{\delta,k}$
are real, and the transform $\varphi_\delta$ is a combination
of exponentials without oscillations.\label{rem:kkompl} Already at $k=3$
and $k=4$, in the example considered, at $\delta=0$ the complex conjugate pairs
$3.82976\pm1.18122\,\mathrm{i}$ and $3.19681\pm1.76237\,\mathrm{i}$
appear, respectively. Therefore the result of
Theorem~\ref{thm:spectrum} on the two positive roots applies specifically to
$k=2$. For $k\ge3$ Rouch\'e's theorem, applied to the equation
$\lambda+\delta-cs=\lambda\mu^k/(\mu+s)^k$ ($s=-x$) in the half-plane
$\operatorname{Re}s\ge0$, gives for $\delta>0$ exactly one root with
$\operatorname{Re}x<0$ and exactly $k$ roots with $\operatorname{Re}x>0$,
among which complex conjugate pairs appear (damped oscillations). The matrix
reduction \eqref{eq:kchain} remains valid.

\section{Properties of the Laplace transform and the expansion in powers of $\delta$}
\label{sec:discount}

The properties of $\varphi_\delta(0)$ as a function of $\delta$ follow
from the definition \eqref{eq:def} itself and do not depend on the method
of computation. Recall that a function $g$ is called completely monotone
on $(0,\infty)$ if it is infinitely differentiable and
$(-1)^n g^{(n)}(\delta)\ge0$ for all $n=0,1,2,\ldots$ and all $\delta>0$.

\begin{lemma}[completely monotone functions]\label{lem:mono}
The function $\delta\mapsto\varphi_\delta(0)$ is infinitely
differentiable and completely monotone on $(0,\infty)$:
\[
(-1)^n\frac{\dd^{n}}{\dd\delta^{n}}\varphi_\delta(0)
=\mathbb{E}_0\bigl[T^{n}e^{-\delta T};\,T<\infty\bigr]\ge0,
\qquad n=0,1,2,\ldots.
\]
In particular, $\varphi_\delta(0)$ decreases, is convex, and lies below
its value at $\delta=0$.
\end{lemma}
\begin{proof}
By \eqref{eq:def} we have
$\varphi_\delta(0)=\int_{[0,\infty)}e^{-\delta t}\nu(\dd t)$,
where $\nu$ is the positive measure of the distribution of the time of
ruin on the set $\{T<\infty\}$. That the Laplace transform of a positive
measure is a completely monotone function is checked by differentiating
under the integral sign (the converse is Bernstein's theorem,
\cite[Ch.~XIII, \S~4]{Feller}). For a fixed
$\delta>0$ pick $0<\delta_0<\delta$. Then
$t^ne^{-\delta t}\le C_{n,\delta,\delta_0}e^{-\delta_0 t}$, and the
integral $\int e^{-\delta_0t}\nu(\dd t)=\varphi_{\delta_0}(0)\le1$ is
finite. This justifies the differentiation under the integral sign.
\end{proof}

\begin{corollary}\label{cor:mono_bounds}
In the model case $\lambda=\mu=3$, $c=6$, for all $\delta\ge0$ the
following estimates hold:
\[
0\le\varphi_\delta(0)\le\varphi_0(0)=\frac13,\qquad
\varphi_\delta(0)\ge\varphi_0(0)+\delta\,\partial_\delta\varphi_\delta(0)\big|_{\delta=0+}=\frac13-\frac{\delta}{24}.
\]
The second estimate follows from convexity; for $\delta>8$ it is trivial.
\end{corollary}

We turn to the analyticity at zero. Put
$d_j(\delta)=cx_j(\delta)+\lambda+\delta$. The substitution of the
two-exponential representation \eqref{eq:twomode} into the boundary system
\eqref{eq:bc} gives
\begin{equation}\label{eq:coeffdelta}
A_1^{(\delta)}=\frac{\lambda x_2}{(x_2-x_1)d_1},\qquad
A_2^{(\delta)}=-\frac{\lambda x_1}{(x_2-x_1)d_2},
\end{equation}
and hence
\begin{equation}\label{eq:phi0delta}
\varphi_\delta(0)=
\frac{\lambda\bigl(c(x_1+x_2)+\lambda+\delta\bigr)}{d_1d_2}.
\end{equation}

\begin{proposition}[analyticity of the two-exponential representation]
\label{prop:analytic}
The branches $x_1(\delta),x_2(\delta)$, the coefficients $A_{1,2}^{(\delta)}$
and the function $\varphi_\delta(u)$, for each fixed $u$, are analytic in
$\delta$ in a neighborhood of zero.
\end{proposition}
\begin{proof}
The positive roots $x_1(0)$ and $x_2(0)$ are simple and distinct by
Theorem~\ref{thm:spectrum}; therefore the implicit function theorem in the
complex domain gives analytic branches $x_1(\delta)$ and $x_2(\delta)$,
which remain distinct in a small neighborhood of zero.
Moreover, $d_j\ne0$ on the corresponding root: otherwise, from
$q_\delta(x_j)=(cx_j+\lambda+\delta)(\mu-x_j)^2-\lambda\mu^2$ it would
follow that $q_\delta(x_j)=-\lambda\mu^2\ne0$. Therefore the denominators
in \eqref{eq:coeffdelta} and \eqref{eq:phi0delta} do not vanish, which
proves the statement.
\end{proof}

For negative $\delta$ inside the domain of analyticity, the formulas
\eqref{eq:twomode}--\eqref{eq:phi0delta} define the continuation of the
branches selected at $\delta\ge0$. For $\delta<0$ boundedness no longer
excludes the third exponential term; therefore the identification of this
continuation with
the expectation \eqref{eq:def} requires a separate justification; at
$u=0$ it is given after Proposition~\ref{prop:radius} (at $u>0$ the
argument is the same).

Theorem~\ref{thm:hier} and Proposition~\ref{prop:analytic} give the moment
series
\[
\varphi_\delta(0)=\sum_{n\ge0}\frac{(-\delta)^{n}}{n!}\,
\mathbb{E}_0\bigl[T^{n};T<\infty\bigr],
\]
and in particular, in the model case,
\begin{equation}\label{eq:dseries}
\varphi_\delta(0)=\frac13-\frac{\delta}{24}+\frac{25\,\delta^{2}}{3456}
-\frac{\delta^{3}}{6}\,\mathbb{E}_0[T^{3};T<\infty]+O(\delta^{4}).
\end{equation}

\begin{proposition}\label{prop:radius}
In the model case the singularity of the constructed analytic continuation
of the function $\delta\mapsto\varphi_\delta(0)$ closest to zero lies at
$\delta=-\gamma_1$, where
\[
\begin{gathered}
\gamma_1=c(\mu-\rho_3)+\lambda-\frac{\lambda\mu^{2}}{\rho_3^{2}},\\
\rho_3=\Bigl(\frac{2\lambda\mu^{2}}{c}\Bigr)^{1/3},\qquad
x_*=\mu-\rho_3 .
\end{gathered}
\]
Numerically $\gamma_1=21-9\sqrt[3]{9}=2.2792455925\ldots$ and
$x_*=0.9199161769\ldots$.
\end{proposition}
\begin{proof}
A singularity of the formulas
\eqref{eq:coeffdelta}--\eqref{eq:phi0delta} can arise only from a
coalescence of branches, since $d_j$ does not vanish at the root. Hence the
equations $q_\delta(x)=\partial_xq_\delta(x)=0$ must hold simultaneously. From
$q_\delta(x)=(cx+\lambda+\delta)(\mu-x)^{2}-\lambda\mu^{2}$ and
$\partial_xq_\delta=c(\mu-x)^{2}-2(cx+\lambda+\delta)(\mu-x)$, at
$y=\mu-x\ne0$ we obtain $cy^3=2\lambda\mu^2$ and
$\delta=\lambda\mu^2/y^2-c(\mu-y)-\lambda$. The real solution
$y=\rho_3$ gives $x_*=\mu-\rho_3$ and $\delta=-\gamma_1$. The case $y=0$
is impossible, since $q_\delta(\mu)=-\lambda\mu^2\ne0$. The two remaining
complex conjugate solutions are obtained by replacing $\rho_3$ with the other
roots of the cubic equation; the modulus of the corresponding $\delta$
equals $34.41805468\ldots>\gamma_1$. As we move from zero to the left, at
the real singularity the branch $x_1$ coalesces with the branch continuing
the zero root, while the branch $x_2$ remains simple.

The point found is a genuine singularity, not a removable one. There
$q_{xx}=-3c\rho_3\ne0$ and
$\partial_\delta q=(\mu-x_*)^2\ne0$, hence
\[
x_1(\delta)=x_*+a(\delta+\gamma_1)^{1/2}
+O(\delta+\gamma_1),\qquad a\ne0.
\]
If $F$ is the right-hand side of \eqref{eq:phi0delta}, then direct
differentiation gives at the critical point
$\partial F/\partial x_1=-\lambda c^2x_2/(d_1^2d_2)\ne0$.
Hence the branch point of square-root type survives.
\end{proof}

For illustration, the values of the analytic continuation itself are
$\varphi_{-0.5}(0)=0.356219$, $\varphi_{-1}(0)=0.384591$,
$\varphi_{-1.5}(0)=0.422364$, $\varphi_{-2}(0)=0.482496$ and
$\varphi_{-2.25}(0)=0.551167$ to the digits shown.
For $-\gamma_1<\delta\le0$ this continuation coincides with
$\mathbb E_0[e^{|\delta|T};T<\infty]$: by Theorem~\ref{thm:hier}
the coefficients of the moment series equal
$\mathbb E_0[T^n;T<\infty]/n!\ge0$, the series converges in the disk
$|\delta|<\gamma_1$, and it remains to apply Tonelli's theorem. For
$|\delta|>\gamma_1$ the corresponding exponential moment is infinite, so
that $\gamma_1$ is the rate of exponential decay
of the distribution of $T$ on $\{T<\infty\}$.

\begin{table}[ht]
\centering
\caption{Extraction of the third moment: the quantities
$Q(\delta)=6\bigl(\textstyle\frac13-\frac{\delta}{24}+\frac{25\delta^{2}}{3456}-\varphi_\delta(0)\bigr)\big/\delta^{3}$
approach $\mathbb{E}_0[T^3;T<\infty]$ linearly in $\delta$.}
\label{tab:m3}
\begin{tabular}{|c|c|}
\hline
$\delta$ & $Q(\delta)$\\ \hline
$0.005$ & $0.0100411$\\
$0.01$ & $0.0100272$\\
$0.02$ & $0.0099995$\\
$0.03$ & $0.0099721$\\
$0.05$ & $0.0099177$\\ \hline
\end{tabular}
\end{table}

\begin{proposition}\label{prop:m3}
In the model case, exactly,
\[
\begin{aligned}
\mathbb{E}_0\bigl[T^{3};T<\infty\bigr]
&=\frac{139}{13824}=0.01005497685185\ldots,\\
\mathbb{E}_0\bigl[T^{3}\mid T<\infty\bigr]&=\frac{139}{4608}.
\end{aligned}
\]
\end{proposition}
\begin{proof}
We expand both positive roots $x_1(\delta),x_2(\delta)$ from
\eqref{eq:char} to the third order in $\delta$ and substitute the two
decaying exponential terms from \eqref{eq:twomode} into the boundary
system \eqref{eq:bc}.
Solving for $\varphi_\delta(0)=A_1^{(\delta)}+A_2^{(\delta)}$ gives
\[
\varphi_\delta(0)=\frac13-\frac{\delta}{24}
+\frac{25\delta^2}{3456}-\frac{139\delta^3}{82944}+O(\delta^4).
\]
The comparison with \eqref{eq:dseries} shows that
$\mathbb E_0[T^3;T<\infty]=139/13824$. The conditional moment is obtained
by dividing by $\mathbb P_0(T<\infty)=1/3$. The numerical values of
Table~\ref{tab:m3} agree with this exact result.
\end{proof}

\medskip
\noindent\textbf{Remark 5.} The approximation by the terms up to the second
order in \eqref{eq:dseries} has a cubic remainder.\label{rem:seriescheck}
Indeed, the difference between $\varphi_\delta(0)$ and
$1/3-\delta/24+25\delta^{2}/3456$ at $\delta=0.01,0.05,0.1,0.2$
equals respectively $-1.67\cdot10^{-9}$, $-2.07\cdot10^{-7}$,
$-1.63\cdot10^{-6}$, $-1.27\cdot10^{-5}$ to the digits
shown. The division of these values by $-\delta^3$ gives quantities close
to $\mathbb{E}_0[T^3;T<\infty]/6=1.68\cdot10^{-3}$.

\section{Conclusion and open questions}\label{sec:open}

For claims with the Erlang distribution of shape parameter~2, we obtain a
closed description of the Laplace transform of the time of ruin: the two
admissible exponential terms are determined by the characteristic equation
and the boundary conditions. This makes it possible to compute both the ruin probability and the
moments of the time of ruin, keeping the full truncated convolution in all
the checks.

For the general $E_k(\mu)$ distribution a finite-dimensional system is
constructed; the structure of the roots, however, is more complicated: in
addition to real exponential terms, complex conjugate pairs are possible, and the
explicit formulas for the coefficients of the $k$ decaying terms require a
separate analysis. Extensions to models with proportional investment and
two-sided jumps, where the variable coefficients or additional unbounded
terms require a different boundary value setting, are separate directions
of further work.

\section*{Declarations}

\paragraph{Conflict of interest.} The author declares no conflict of interest.

\paragraph{Availability of data and materials.} No datasets were generated or analyzed during the current study.

\end{document}